# HITTING PROPERTIES OF PARABOLIC S.P.D.E.'S WITH REFLECTION

By Robert C. Dalang[1]  C. Mueller[2]  and  L. Zambotti

*École Polytechnique Fédérale, University of Rochester and Politecnico di Milano*

We study the hitting properties of the solutions $u$ of a class of parabolic stochastic partial differential equations with singular drifts that prevent $u$ from becoming negative. The drifts can be a reflecting term or a nonlinearity $cu^{-3}$, with $c > 0$. We prove that almost surely, for all time $t > 0$, the solution $u_t$ hits the level 0 only at a finite number of space points, which depends explicitly on $c$. In particular, this number of hits never exceeds 4 and if $c > 15/8$, then level 0 is not hit.

**1. Introduction.** We consider the nonnegative solutions $u$ of the following classes of stochastic partial differential equations (s.p.d.e.'s) driven by space–time white noise: the heat equation with repulsion from 0 (introduced in [18])

$$(1.1) \qquad \frac{\partial u}{\partial t} = \frac{1}{2}\frac{\partial^2 u}{\partial x^2} + \frac{c}{u^3} + \frac{\partial^2 W}{\partial t\,\partial x}, \qquad c > 0;$$

and the heat equation with reflection at 0 (introduced in [11])

$$(1.2) \qquad \frac{\partial u}{\partial t} = \frac{1}{2}\frac{\partial^2 u}{\partial x^2} + \frac{\partial^2 W}{\partial t\,\partial x} + \eta, \qquad c = 0.$$

In (1.1) and (1.2), $x \in [0, 1]$, $W = (W(t, x), t \geq 0, x \in \mathbb{R})$ is a Brownian sheet on a complete probability space $(\Omega, \mathcal{F}, \mathbb{P})$, the continuous solution $u = (u_t(x), t \geq 0, x \in [0,1])$ satisfies homogeneous Dirichlet boundary conditions at 0 and 1, and $u \geq 0$ on $[0, \infty) \times [0, 1]$. Moreover, in (1.2), $\eta$ is a nonnegative

Received August 2004; revised April 2005.
[1]Supported in part by the Swiss National Foundation for Scientific Research.
[2]Supported in part by an NSF grant.

measure on $(0,\infty) \times (0,1)$ that is supported by $\{(t,x): u_t(x) = 0\}$ and is called the reflecting measure.

The aim of this paper is to study the random contact set

$$\mathcal{Z} = \mathcal{Z}^{(c)} := \{(t,x): u_t(x) = 0\}$$

for solutions $u = u^{(c)}$ of (1.1) for $c > 0$, respectively, (1.2) for $c = 0$. Notice that (1.1) and (1.2) are stochastic obstacle problems: indeed $u \geq 0$ solves an s.p.d.e. outside the contact set $\mathcal{Z}$, which itself is determined by $u$. Since the drifts of our equations become singular as $u$ approaches 0, then we expect $\mathcal{Z}$ to be smaller than level sets $\{(t,x): u_t(x) = a\}$ with $a > 0$.

An important property of $(u^{(c)})_{c \geq 0}$ is monotonicity in $c$ (see the proof of Lemma 3.1): for given initial and boundary conditions, if $c \geq c' \geq 0$, then a.s. $u^{(c)} \geq u^{(c')}$, so that a.s. $\mathcal{Z}^{(c)} \subseteq \mathcal{Z}^{(c')}$. Therefore, it is natural to conjecture that there exists a $c_0 > 0$, possibly random, such that $\mathcal{Z}^{(c)} = \varnothing$ for all $c > c_0$. On the other hand, it is not easy to guess the behavior of $\mathcal{Z}^{(c)}$ for small $c$.

In this paper we study the cardinality of the $x$ sections of the random set $\mathcal{Z}$: that is, for all $t > 0$ we define

$$\zeta_t := |\{x \in (0,1): u_t(x) = 0\}|,$$

where $|\cdot|$ denotes the cardinality of a set. Then we consider the random variable

$$(1.3) \qquad \zeta = \zeta^{(c)} := \sup_{t \in (0,1]} \zeta_t.$$

Notice that in the definition of $\zeta_t$, we exclude $x \in \{0,1\}$, because there $u = 0$ by the boundary conditions. By the monotonicity in $c$, we have $\zeta^{(c)} \leq \zeta^{(c')}$ if $c \geq c'$.

Our main results give much more precise information about $\zeta$. First we prove that for all $c \geq 0$, $\zeta^{(c)} \leq 4$ a.s. This is rather surprising, due to the wild oscillations of the space–time white noise and to the zero boundary conditions.

Notice that for all $c \geq 0$, there exists a unique $\delta \geq 3$ such that

$$(1.4) \qquad c = c_\delta := \frac{(\delta-3)(\delta-1)}{8}.$$

In [17] and [18], it is proven that the process $(u_t)_{t \geq 0}$ is stationary if and only if $u_0$ is distributed like a Bessel bridge of dimension $\delta$; see [14], Chapter XI. In particular, we have at hand the explicit law of $u_t(x)$ for stationary $(u_t)_{t \geq 0}$, which turns out to be crucial for obtaining our results.

The second result of this paper states that for all $\delta > 6$, $\zeta = 0$ a.s. Therefore, if $c > 15/8$, then the contact set is empty.

As $\delta$ decreases from 6 to 3, we have the intermediate behavior

$$(1.5) \qquad \zeta(\delta) \leq \frac{4}{\delta - 2} \qquad \text{a.s.}$$

In particular, since $\zeta$ is an integer, $\zeta = 0$ a.s. for $\delta > 6$, $\zeta \leq 1$ a.s. for all $\delta \in (4, 6]$, $\zeta \leq 2$ a.s. for all $\delta \in (10/3, 4]$, $\zeta \leq 3$ a.s. for all $\delta \in (3, 10/3]$ and $\zeta \leq 4$ a.s. for $\delta = 3$.

We also give nontrivial lower bounds for $\zeta$. Indeed, we prove that with positive probability, $\zeta \geq 1$ for all $\delta \in [3, 5]$ and $\zeta \geq 3$ for $\delta = 3$. The latter result is particularly interesting, for the following reason. In [17], it was proved that for $\delta = 3$ or, equivalently, for $c = 0$, almost surely,

$$(1.6) \qquad \zeta_t = 1 \qquad \text{for } \eta(dt \times (0,1))\text{-a.e. } t.$$

Therefore, generically $\zeta_t = 1$ at typical times. By the result of this paper, (1.6) is not optimal and the set of times $t > 0$ such that $\zeta_t > 1$ is a.s. negligible for $\eta$, but nonempty with positive probability.

We recall that Mueller [7] and Mueller and Pardoux [8] considered the s.p.d.e. with periodic boundary conditions,

$$\begin{cases} \dfrac{\partial \hat{u}}{\partial t} = \dfrac{\partial^2 \hat{u}}{\partial \theta^2} + \hat{u}^{-\alpha} + g(\hat{u}) \dfrac{\partial^2 W}{\partial t \, \partial \theta}, & t \geq 0, \qquad \theta \in \mathbb{S}^1 := \mathbb{R}/\mathbb{Z}, \\ \hat{u}(0, \cdot) = \hat{u}_0(\cdot), \end{cases}$$

where $\alpha > 0$, $\hat{u}_0 : \mathbb{S}^1 \mapsto \mathbb{R}$ is continuous, $\inf \hat{u}_0 > 0$ and $g$ satisfies suitable growth conditions, and proved that $\alpha = 3$ is the critical exponent for $\hat{u}$ to hit zero in finite time. More precisely, the following statements were proved.

1. If $\alpha > 3$, then a.s. $\hat{u}(t, \theta) > 0$ for all $t \geq 0$ and $\theta \in \mathbb{S}^1$.
2. If $\alpha < 3$, then with positive probability, there exist $t > 0$ and $\theta \in \mathbb{S}^1$ such that $\hat{u}(t, \theta) = 0$.

Existence for all time of a solution for $\alpha = 3$ was first proved in [18]. In this paper, we prove that, in the critical case $\alpha = 3$, the hitting properties of the solution depend on the constant $c$. This is reminiscent of the behavior of the Bessel processes $(X_t)_{t \geq 0}$, solution of

$$dX_t = \frac{\delta - 1}{2 X_t} dt + dB, \qquad X_0 \geq 0,$$

where $\delta > 1$ and $B$ is a standard Brownian motion. Indeed, it is well known that $X$ hits 0 with positive probability if and only if $\delta < 2$; see [14], Chapter XI, Section 1.

Further questions addressed in this paper concern the study of similar hitting properties for multidimensional solutions of linear s.p.d.e.'s, which continues the work of Mueller and Tribe [9]. For this class of Gaussian processes, we derive optimal results.

**2. Main results.** We define $C^+ := \{\bar{u}:[0,1] \mapsto [0,\infty): \bar{u}$ is continuous, $\bar{u}(0) = \bar{u}(1) = 0\}$ and we consider a Brownian sheet $(W(t,x): t \geq 0, x \in [0,1])$ and the associated filtration $(\mathcal{F}_t, t \geq 0)$, where

$$\mathcal{F}_t = \sigma\{W(s,x), s \leq t, x \in [0,1]\} \vee \mathcal{N}$$

and $\mathcal{N}$ is the $\sigma$-field generated by all $P$-null sets. For any $\delta > 0$ and $\bar{u} \in C^+$, we consider the unique continuous nonnegative solutions $(u_t(x): t \geq 0, x \in [0,1])$ of the following s.p.d.e.'s:

(2.1) $\quad \delta \in (3, \infty), \quad \begin{cases} \dfrac{\partial u}{\partial t} = \dfrac{1}{2}\dfrac{\partial^2 u}{\partial x^2} + \dfrac{(\delta-1)(\delta-3)}{8u^3} + \dfrac{\partial^2 W}{\partial t\, \partial x}, \\ u_t(0) = u_t(1) = 0, \quad t \geq 0, \\ u_0(x) = \bar{u}(x), \quad x \in [0,1], \end{cases}$

(2.2) $\quad \delta = 3, \quad \begin{cases} \dfrac{\partial u}{\partial t} = \dfrac{1}{2}\dfrac{\partial^2 u}{\partial x^2} + \dfrac{\partial^2 W}{\partial t\, \partial x} + \eta(t,x), \\ u_0(x) = \bar{u}(x), \quad u_t(0) = u_t(1) = 0, \\ u \geq 0, \quad d\eta \geq 0, \quad \int u\, d\eta = 0. \end{cases}$

Rigorous meanings of the equations (2.1) and (2.2) are, respectively, given in [18] and [11], where existence and uniqueness of solutions are also proved. We recall that in (2.1), the unique solution satisfies $u^{-3} \in L^1_{\mathrm{loc}}((0,\infty) \times (0,1))$. Moreover, in (2.2), the nonnegative measure $\eta$ is a reflecting term, with support included in $\{(t,x): u_t(x) = 0\}$. In [18], it is proved that the solution of (2.1) converges a.s. to the solution of (2.2) as $\delta \searrow 3$. For this reason, we interpret (2.2) as the case $\delta = 3$ of (2.1).

The main results of this paper are the following two theorems.

THEOREM 2.1. *For all $\bar{u} \in C^+$, the following statements hold.*

(a) *For $\delta > 6$, the probability that there exist $t > 0$ and $x \in (0,1)$ such that $u_t(x) = 0$ is 0.*

(b) *For $\delta > 4$, the probability that there exist $t > 0$ and $\{x_i, i=1,2\} \subset (0,1)$, $x_1 < x_2$, such that $u_t(x_i) = 0$, $i = 1, 2$, is 0.*

(c) *For $\delta > \frac{10}{3}$, the probability that there exist $t > 0$ and $\{x_i, i=1,2,3\} \subset (0,1)$, $x_1 < x_2 < x_3$, such that $u_t(x_i) = 0$, $i = 1, 2, 3$, is 0.*

(d) *For $\delta > 3$, the probability that there exist $t > 0$ and $\{x_i, i=1,\ldots,4\} \subset (0,1)$, $x_1 < \cdots < x_4$, such that $u_t(x_i) = 0$, $i = 1,\ldots,4$, is 0.*

(e) *For $\delta = 3$, the probability that there exist $t > 0$ and $\{x_i, i=1,\ldots,5\} \subset (0,1)$, $x_1 < \cdots < x_5$, such that $u_t(x_i) = 0$, $i = 1,\ldots,5$, is 0.*

THEOREM 2.2. *For all $\bar{u} \in C^+$, the following statements hold.*

(a) *For all $\delta \in [3,5]$, with positive probability there exist $t > 0$ and $x \in (0,1)$ such that $u_t(x) = 0$.*

(b) *For $\delta = 3$, with positive probability there exist $t > 0$ and $\{x_1, x_2, x_3\} \subset (0,1)$, $x_1 < x_2 < x_3$, such that $u_t(x_i) = 0$, $i = 1, 2, 3$.*

Notice that these results are optimal only for $\delta \in (4, 5]$, since it is only for such $\delta$ that they imply that the upper bound for $\zeta$ [defined in (1.3)] is attained with positive probability.

We recall that Mueller and Tribe [9] have defined the *stationary pinned string*, that is, the solution $U_t(x) \in \mathbb{R}^d$, $d \in \mathbb{N}$, of

$$(2.3) \qquad \frac{\partial U_t}{\partial t} = \frac{1}{2}\frac{\partial^2 U_t}{\partial x^2} + \frac{\partial^2 W_d}{\partial t\, \partial x}, \qquad t > 0, x \in \mathbb{R},$$

where $W_d = (W^1, \ldots, W^d)$, $\{W^i\}_{i=1,\ldots,d}$ is an independent sequence of copies of $W$ and $(U_0(x): x \in \mathbb{R})$ is a two-sided $\mathbb{R}^d$-valued Brownian motion independent of $W_d$ and satisfying

$$U_0(0) = 0, \qquad \mathbb{E}[(U_0(x) - U_0(y))^2] = |x - y|^2.$$

In particular,

$$(2.4) \qquad U_t(x) = \int_{\mathbb{R}} G_t(x - z) U_0(z)\, dz + \int_0^t \int_{\mathbb{R}} G_s(x - z)\, W_d(ds, dz),$$

where $G_t$ is the density of the Gaussian distribution with mean zero and variance $t$. The following result identifies the dimensions in which the stationary pinned string hits points.

THEOREM 2.3 ([9], Theorem 1). *The probability that there exist $t > 0$ and $x \in \mathbb{R}$ such that $U_t(x) = 0$ is positive if and only if $d \leq 5$.*

In this paper we complete this result as follows. First, following the definition of $\zeta$ in (1.3), we introduce the random variable

$$(2.5) \qquad Z = Z(d) := \sup_{t \in (0,1]} |\{x \in \mathbb{R} : U_t(x) = 0\}|,$$

where $|\cdot|$ again denotes cardinality.

THEOREM 2.4.

1. *For $d \geq 4$, the probability that there exist $t > 0$ and $x_1 < x_2$ such that $U_t(x_i) = 0$, $i = 1, 2$, is 0. [In fact, for $d \geq 4$, $Z(d) \leq 1$ a.s. and by Theorem 2.3, $P\{Z(d) = 1\} > 0$ if $d \in \{4, 5\}$ and $Z(d) = 0$ a.s. if $d \geq 6$.]*
2. *The probability that there exist $t > 0$ and $x_1 < x_2 < x_3$ such that $U_t(x_i) = 0$, $i = 1, 2, 3$, is positive if and only if $d \leq 3$. In addition, $Z(3) \leq 3$ a.s.*
3. *If $d = 2$, then for all $k \in \mathbb{N}$, with positive probability, there exist $t > 0$ and $x_1 < \cdots < x_k$ such that $U_t(x_i) = 0$, $i = 1, \ldots, k$.*

Notice that for the Gaussian process $U$ and for $d \geq 3$, our upper bounds are attained with positive probability. Notice also that Theorems 2.2 and 2.4 are related: in fact, for $d, k \in \mathbb{N}$ with $d \geq 3$, the following implications hold:

$$\mathbb{P}(Z(d) = k) > 0 \quad \Longrightarrow \quad \mathbb{P}(\zeta(\delta) = k) > 0 \qquad \forall \delta \leq d,$$
$$\mathbb{P}(Z(d) = k) = 0 \quad \Longrightarrow \quad \mathbb{P}(\zeta(\delta) = k) = 0 \qquad \forall \delta > d.$$

These relationships can be explained with a result of [17] for $c = 0$ and [18] for $c > 0$, relating (1.1), (1.2) and (2.3) for $\delta = d \in \mathbb{N}$: see the proofs of Theorems 2.1 and 2.2 below.

Although our approach does not yield optimal results for the nonlinear equations (1.1) and (1.2), on the basis of Theorem 2.4 we can propose the following conjectures:

1. We conjecture that $\zeta(\delta) < \frac{4}{\delta - 2}$ a.s. for all $\delta \geq 3$.
2. We conjecture that $\mathbb{P}(\zeta(\delta) = 1) > 0$ for $\delta \in [4, 6)$, $\mathbb{P}(\zeta(\delta) = 2) > 0$ for $\delta \in [10/3, 4)$ and $\mathbb{P}(\zeta(\delta) = 3) > 0$ for $\delta \in [3, 10/3)$.

Part 1 would improve (1.5) and Theorem 2.1 when $\delta \in \{3, 10/3, 4, 6\}$, part 2 would improve Theorem 2.2 for $\delta \in (3, 4) \cup (5, 6)$ and these bounds would be optimal.

This paper is organized as follows. In Section 3 we study the same s.p.d.e.'s as (2.1) and (2.2), but with positive boundary conditions for the former. This makes it possible to establish some Hölder continuity properties of the solution and to prove the analog of Theorem 2.1 in this case. In Section 4 we use the results of Section 3 and some comparison theorems to prove Theorem 2.1. In Section 5 we turn to the vector-valued linear equation (2.3), proving Theorem 2.4. In Section 6 we use Theorem 2.4 to establish Theorem 2.2.

**3. Hölder continuity and a variant on Theorem 2.1.** In this section we prove a variant on Theorem 2.1 in which the boundary conditions of the s.p.d.e. (2.1) are positive; those of (2.2) may be positive or may vanish. Our approach to proving this theorem uses a classical discretization technique [4]. After choosing a grid in the rectangle $[0, T] \times [0, 1]$, we perform two steps: first, we prove that the probability of finding a point on the grid where $u$ is close to 0 is small; second, we control the oscillations of $u$, proving that the result on the grid extends to the whole rectangle.

The first step is based on the explicit knowledge of the invariant measure of equations (2.1) and (2.2), obtained in [17] and [18]: indeed, in the stationary case the distribution of $u_t(x)$ is known for fixed $(t, x)$ in the grid.

The second step is based on an estimate of the Hölder regularity of $u$. This issue is nontrivial since the nonlinearities in (2.1) and (2.2) become singular as $u \to 0$. In fact, we can prove that $u$ is Hölder-continuous in space, but as far as time regularity is concerned, only our lower bound is optimal: since the singular terms are positive, $u$ does not *decrease* too quickly. See Lemma 3.1 and, in particular (3.5), as well as Remark 3.7.

Let $\delta \geq 3$, $[b,c] \subseteq [0,1]$ and $a \geq 0$, and denote by $(v_t(x) : t \geq 0, x \in [b,c])$ the unique solution of

$$(3.1) \quad \delta \in (3,\infty), \quad \begin{cases} \dfrac{\partial v}{\partial t} = \dfrac{1}{2} \dfrac{\partial^2 v}{\partial x^2} + \dfrac{(\delta-1)(\delta-3)}{8v^3} + \dfrac{\partial^2 W}{\partial t \, \partial x}, \\ v_t(b) = v_t(c) = a, \quad t \geq 0, \\ v_0(x) = \overline{v}(x), \quad x \in [b,c], \end{cases}$$

$$(3.2) \quad \delta = 3, \quad \begin{cases} \dfrac{\partial v}{\partial t} = \dfrac{1}{2} \dfrac{\partial^2 v}{\partial x^2} + \dfrac{\partial^2 W}{\partial t \, \partial x} + \zeta(t,x), \\ v_0(x) = \overline{v}(x), \quad v_t(b) = v_t(c) = a, \\ v \geq 0, \quad d\zeta \geq 0, \quad \int v \, d\zeta = 0, \end{cases}$$

where $\overline{v} : [b,c] \mapsto \mathbb{R}$ is continuous nonnegative with $\overline{v}(b) = \overline{v}(c) = a$. Clearly, $u = v$ if $a = 0$ and $[b,c] = [0,1]$.

For $b < c$ and $\beta > 0$, let $C^\beta([b,c])$ denote the space of Hölder-continuous functions on $[b,c]$ with Hölder exponent $\beta$, equipped with the norm

$$\|\overline{v}\|_\beta := \sup_{x \in [b,c]} |\overline{v}(x)| + \sup_{b < x < y < c} \frac{|\overline{v}(x) - \overline{v}(y)|}{|x-y|^\beta}.$$

In the proof of Theorem 2.1, the following lemma plays a key role.

LEMMA 3.1. *Let $\delta \geq 3$, $a \geq 0$ satisfy* (I) *or* (II)*, where*

$$(3.3) \qquad \text{(I)} \quad \delta = 3, a \geq 0; \qquad \text{(II)} \quad \delta > 3, a > 0.$$

*Let $(v_t(x) : t \geq 0, x \in [b,c])$ satisfy (3.1) or (3.2). Then for all $\beta \in (0, 1/2)$ and $T > 0$, if $\overline{v} \in C^\beta([b,c])$, then there exists a finite random variable $\gamma_v$ such that*

$$(3.4) \quad |v_t(x) - v_t(y)| \leq \gamma_v |x-y|^\beta, \qquad x, y \in [0,1], \qquad T \geq t \geq 0,$$

*and*

$$(3.5) \quad v_t(x) - v_s(x) \geq -\gamma_v (t-s)^{\beta/2}, \qquad T \geq t \geq s \geq 0, \qquad x \in [0,1].$$

We postpone the proof of Lemma 3.1 to the end of this section. Let $(g_t(x,y) : t > 0, x, y \in [b,c])$ be the Green function of the heat equation with homogeneous Dirichlet boundary conditions

$$\begin{cases} \dfrac{\partial g}{\partial t} = \dfrac{1}{2} \dfrac{\partial^2 g}{\partial x^2}, & t > 0, \quad x \in (b,c), \\ g_t(b,y) = g_t(c,y) = 0, & t > 0, \quad y \in (b,c), \\ g_0(x,y) = \delta_x(y), & x \in (b,c), \end{cases}$$

where $\delta_x$ is the Dirac mass at $x \in (b,c)$.

REMARK 3.2. As proven in [16], for the stochastic convolution,

$$(3.6) \qquad S_t^{(\bar v)}(x) := \int_b^c g_t(x,y)\bar v(y)\,dy + \int_0^t \int_b^c g_{t-s}(x,y)W(ds,dy),$$

if $\bar v \in C^\beta([b,c])$, then there exists a finite random variable $\gamma_S$ such that a.s., for all $t,s \in [0,T]$, $x,y \in [b,c]$,

$$(3.7) \qquad |S_t^{(\bar v)}(x) - S_s^{(\bar v)}(y)| \leq \gamma_S(|t-s|^{\beta/2} + |x-y|^\beta).$$

By (3.4), $v$ satisfies the same Hölder continuity in space as $S_t^{(\bar v)}(\cdot)$. On the other hand, the singular drift $v^{-3}$ might produce worse behavior in time, in particular around $(t,x)$ such that $v_t(x) = 0$. Nevertheless, by (3.5), $t \mapsto v_t(x)$ cannot *decrease* more quickly than $t \mapsto S_t^{(\bar v)}(x)$.

We denote by $(X_\theta : \theta \in [b,c])$ a Bessel bridge of dimension $\delta$ between $a$ and $a$ over the interval $[b,c]$ (see [14]). We shall exploit the relationship of this bridge with the Bessel process $Y^{(\delta)}$ of dimension $\delta$. Let $p_t(x,y)$ denote the transition semigroup of $Y^{(\delta)}$. We recall that for $x > 0$, $y \geq 0$ and $t > 0$,

$$(3.8) \qquad p_t(x,y) := \frac{1}{t}\left(\frac{y}{x}\right)^{(\delta/2)-1} y\exp\left(-\frac{x^2+y^2}{2t}\right)I_{(\delta/2)-1}\left(\frac{xy}{t}\right),$$

where $I$ is the modified Bessel function and for $x=0$,

$$p_t(0,y) = \frac{1}{2^{\delta/2-1}t^{\delta/2}\Gamma(\delta/2)}y^{\delta-1}e^{-y^2/(2t)};$$

see [14], Chapter XI, Section 1. We note for future reference that for $x \geq 0$, $I_\nu(x) = x^\nu \lambda_\nu(x)$ with $\lambda_\nu$ locally bounded and $\lambda_\nu(0) > 0$. In particular, for all $t_0 > 0$, there exists a constant $C$ such that

$$(3.9) \qquad p_t(x,y) \leq Cy^{\delta-1} \qquad \forall t \geq t_0, \qquad x,y \in [0,1].$$

We recall that the laws of $(X_\theta : \theta \in [b,(b+c)/2])$ and $(Y_\theta^{(\delta)} : \theta \in [b,(b+c)/2])$ are mutually absolutely continuous. Indeed, let $b=0$ for simplicity. By the Markov property, for any bounded functional $\Phi$,

$$(3.10) \qquad E(\Phi(X_\theta, \theta \leq c/2)) = E(\Phi(Y_\theta^{(\delta)}, \theta \leq c/2)\tilde p_{c/2}(Y_{c/2}^{(\delta)}, a)),$$

where

$$(3.11) \qquad \tilde p_{c/2}(y,a) = \frac{p_{c/2}(y,a)}{p_c(0,a)} \qquad \text{if } a \neq 0$$

and $\tilde p_{c/2}(y,0) = \lim_{a\downarrow 0}\tilde p_{c/2}(y,a) = \exp(-y^2/c)/(c/2)$.

We now recall the following result, proved in [17] for $\delta = 3$ and in [18] for $\delta > 3$.

PROPOSITION 3.3. *For any $\delta \geq 3$, $v$ is stationary if and only if $(\overline{v}(x) : x \in [b,c])$ is distributed like $X$ and independent of $W$.*

We now prove the following lemma.

LEMMA 3.4. *For all $\delta \geq 3$ and $\beta \in (0, 1/2)$, there exists a finite real random variable $\gamma_X$ such that a.s.*

$$|X_\theta - X_{\theta'}| \leq \gamma_X |\theta - \theta'|^\beta, \qquad \theta, \theta' \in [b, c].$$

PROOF. Without loss of generality, we can suppose that $b = 0$. Let $Y^{(\delta)}$ be a Bessel process of dimension $\delta$ with $Y_0^{(\delta)} = a$. Since the laws of $(X_\theta : \theta \in [0, c/2])$ and $(Y_\theta^{(\delta)} : \theta \in [0, c/2])$ are mutually absolutely continuous and the law of $X$ is invariant under the time reversal $\theta \mapsto c - \theta$, it is enough to prove the Hölder continuity of $Y^{(\delta)}$ on $[0, c/2]$.

For $\delta = 3$, the result follows from the equality in law between $Y^{(3)}$ and the modulus of a Brownian motion of dimension 3. Let $(B_\theta)_{\theta \in [0,1]}$ be a standard Brownian motion. We recall that for all $\delta \geq 3$, we can realize $Y^{(\delta)}$ as the unique strong solution of the stochastic differential equation (s.d.e.)

$$Y_\theta^{(\delta)} = a + \int_0^\theta \frac{\delta - 1}{2 Y_s^{(\delta)}} ds + B_\theta, \qquad \theta \in [0, 1]$$

(see [14], Chapter XI, Section 1, which also gives the s.d.e. for the square of $Y^{(\delta)}$). Moreover, via standard comparison theorems (see, e.g., [14], Chapter IX, Section 3), which apply to the s.d.e. for the square of Bessel processes, the following monotonicity holds: if $\delta \geq \delta'$, then $Y^{(\delta)} \geq Y^{(\delta')}$ a.s. Now for any $\delta > 3$ and $\theta \leq \theta'$,

$$|Y_{\theta'}^{(\delta)} - Y_\theta^{(\delta)}| \leq |(Y_{\theta'}^{(\delta)} - Y_\theta^{(\delta)}) - (B_{\theta'} - B_\theta)| + |B_{\theta'} - B_\theta|$$

and the first term on the right-hand side is equal to

$$\int_\theta^{\theta'} \frac{\delta - 1}{2 Y_s^{(\delta)}} ds \leq \int_\theta^{\theta'} \frac{\delta - 1}{2 Y_s^{(3)}} ds$$
$$= \frac{\delta - 1}{2} [(Y_{\theta'}^{(3)} - Y_\theta^{(3)}) - (B_{\theta'} - B_\theta)];$$

hence, the result follows from the Hölder continuity of $B$ and $Y^{(3)}$. □

THEOREM 3.5. *Let $\delta$ and $a$ satisfy* (I) *or* (II) *in* (3.3). *If $k \in \mathbb{N}$ satisfies*

(3.12) $$k > \frac{4}{\delta - 2},$$

*then the probability that there exist $t > 0$ and $x_1, \ldots, x_k \in [b, c]$ such that $b < x_1 < \cdots < x_k < c$ and $v_t(x_i) = 0$ for all $i = 1, \ldots, k$, is 0.*

PROOF. First, we notice that it is enough to consider the case of stationary $v$, that is, by Proposition 3.3, to consider $\overline{v}$ distributed like $X$ and independent of $W$. Indeed, for all $n \in \mathbb{N}$, the law of $(v_t : t \geq 1/n)$ for any $\overline{v} \in C^+$ is absolutely continuous with respect to the law of $(v_t : t \geq 1/n)$ with $v$ stationary, since, as proven in [18], page 341, for any $\overline{v} \in C^+$ the law of $v_{1/n} \in C^+$ is absolutely continuous with respect to the law of $X$.

Now, by Lemma 3.4, $\overline{v} \in C^\beta([0,1])$ a.s. for all $\beta \in (0, 1/2)$ and, by Lemma 3.1, $v$ satisfies (3.4) and (3.5).

Let $\mathbb{Q}$ denote the set of rational numbers. For all $\{q_i : i = 1, \ldots, 2k\} \subset \mathbb{Q}$ such that $b < q_1 < \cdots < q_{2k} < c$, we define $Q := [0,1] \times \prod_{i=1}^{k}[q_{2i-1}, q_{2i}]$ and the random set
$$A := \{(t, x_1, \ldots, x_k) \in Q : v_t(x_i) = 0, i = 1, \ldots, k\}.$$

Then the claim will follow if we prove that $\mathbb{P}(A \neq \varnothing) = 0$ for all such $(q_i)_i$.

By (3.12), we can fix $\alpha \in (0,1)$ such that

(3.13) $$4 + 2k - \alpha \delta k < 0.$$

For such $\alpha$, we define the random set
$$A_n := \{(t, x_1, \ldots, x_k) \in Q : v_t(x_i) \leq 2^{-\alpha n}, i = 1, \ldots, k\}.$$

For all $n \in \mathbb{N}$, let
$$G_n := \{(j 2^{-4n}, i_1 2^{-2n}, \ldots, i_k 2^{-2n}) : j, i_1, \ldots, i_k \in \mathbb{Z}\}$$

and consider the events
$$\mathcal{K}_n := \{A_n \cap G_n \neq \varnothing\} \quad \text{and} \quad \mathcal{L}_n := \{A \neq \varnothing, A_n \cap G_n = \varnothing\}.$$

Since $A \subset A_n$ a.s.,
$$\{A \neq \varnothing\} \subseteq \mathcal{K}_n \cup \mathcal{L}_n.$$

To prove that $P\{A \neq \varnothing\} = 0$, we will show that the probabilities of $\mathcal{K}_n$ and $\mathcal{L}_n$ tend to 0 as $n \to \infty$.

*Step* 1. By definition, on $\mathcal{L}_n$ there exists a random $(t, x) \in [0, 1] \times (b, c)$ such that $v_t(x) = 0$ but $A_n \cap G_n = \varnothing$. In particular, on $\mathcal{L}_n$ there exists a random $(s, y) \in \{(j 2^{-4n}, i 2^{-2n}) : j = 1, \ldots, 2^{4n}, i = 1, \ldots, 2^{2n}\}$ such that
$$v_s(y) > 2^{-\alpha n}, \qquad 0 < t - s \leq 2^{-4n}, \qquad |x - y| \leq 2^{-2n}.$$

Let $\beta \in (\alpha/2, 1/2)$. Then on $\mathcal{L}_n$, by (3.4), (3.5) and because $s < t$,
$$2^{-\alpha n} < v_s(y) = v_s(y) - v_t(x)$$
$$= [v_s(y) - v_t(y)] + [v_t(y) - v_t(x)]$$
$$\leq \gamma_v((t-s)^{\beta/2} + |y-x|^\beta) \leq \gamma_v 2^{-2\beta n+1}.$$

Therefore,

$$\mathbb{P}(\mathcal{L}_n) \leq \mathbb{P}(2^{-\alpha n} < \gamma_v 2^{-2\beta n+1}) = \mathbb{P}(\gamma_v > 2^{(2\beta-\alpha)n-1}) \to 0$$

as $n \to \infty$, since $2\beta > \alpha$ and $\gamma_v$ is a.s. finite.

*Step* 2. We set $I_n := G_n \cap Q$. Then, by definition,

$$\mathbb{P}(\mathcal{K}_n) = \mathbb{P}(\exists\, (t, x_1, \ldots, x_k) \in I_n : v_t(x_i) \leq 2^{-\alpha n}, i = 1, \ldots, k).$$

Let $J_n := \{(x_1, \ldots, x_k) : (0, x_1, \ldots, x_k) \in I_n\}$. Then

$$\begin{aligned}
(3.14) \quad \mathbb{P}(\mathcal{K}_n) &\leq \sum_{j=1}^{2^{4n}} \sum_{(x_1,\ldots,x_k) \in J_n} \mathbb{P}(v_{j2^{-4n}}(x_i) \leq 2^{-\alpha n}, i = 1, \ldots, k) \\
&= 2^{4n} \sum_{(\theta_1,\ldots,\theta_k) \in J_n} \mathbb{P}(X_{\theta_i} \leq 2^{-\alpha n}, i = 1, \ldots, k),
\end{aligned}$$

since we have chosen $u$ to be stationary and therefore, for any $t \geq 0$, $v_t$ is distributed like $X$. By (3.10), for $\epsilon > 0$,

$$\mathbb{P}(X_{\theta_i} \leq \epsilon, i = 1, \ldots, k)$$
$$= \int_{[0,\epsilon)^k} \left[ \prod_{i=1}^{k} p_{\theta_i - \theta_{i-1}}(x_{i-1}, x_i) \right] \tilde{p}_{c-\theta_k}(x_k, a) \, dx_1 \cdots dx_k,$$

where $\theta_0 := b$, $x_0 := a$ and $\tilde{p}_{c-\theta_k}(x_k, a)$ is defined in (3.11). We recall that $\theta_i \in [q_{2i-1}, q_{2i}]$, $i = 1, \ldots, k$, and $0 < q_1 < \cdots < q_{2k} < 1$. In all cases, the factor $\tilde{p}_{c-\theta_k}(x_k, a)$ is bounded above and, therefore, by (3.9), there exists a constant $C > 0$ such that for all $(\theta_i)_{i=1,\ldots,k} \in J_n$,

$$\mathbb{P}(X_{\theta_i} \leq \epsilon, i = 1, \ldots, k) \leq C \left[ \int_0^\epsilon x^{\delta-1} \, dx \right]^k \leq C \epsilon^{\delta k}, \qquad \epsilon > 0.$$

Therefore, by (3.14), since the number of elements of $J_n$ is not more than $2^{2kn}$,

$$\mathbb{P}(\mathcal{K}_n) \leq C 2^{4n} 2^{2kn} (2^{-\alpha n})^{\delta k} = C 2^{(4+2k-\alpha\delta k)n} \longrightarrow 0$$

as $n \to \infty$, by (3.13) above. $\square$

In the proof of Lemma 3.1 we need the following result, which is essentially a version of the maximum principle. For $T > 0$ we set $O_T := [0, T] \times [b, c]$ and

$$\|F\|_T := \sup_{O_T} |F|, \qquad F \in C(O_T).$$

LEMMA 3.6. *Let $V \in C^{1,2}(O_T)$ and $\psi, F \in C(O_T)$ with $\psi \leq 0$. Suppose that $V$ solves the equation*

$$\begin{cases} \dfrac{\partial V}{\partial t} = \dfrac{1}{2}\dfrac{\partial^2 V}{\partial x^2} + \psi \cdot V + \psi \cdot F, \\ V_0(x) = 0 \end{cases} \tag{3.15}$$

*with homogeneous Dirichlet or Neumann boundary conditions. Then the following estimate holds:*

$$\|V\|_T \leq \|F\|_T. \tag{3.16}$$

PROOF. We consider first the case of homogeneous Neumann boundary conditions. We denote by $\mathbf{E}_x$ the law of the reflecting Brownian motion $(x_\tau, \tau \geq 0)$ with values in $[b,c]$ started at $x_0 = x \in [b,c]$:

$$x_\tau = x + B_\tau + \tfrac{1}{2}L_\tau^b - \tfrac{1}{2}L_\tau^c, \qquad \tau \geq 0,$$

where $L^\alpha$ is the local time process of $(x_\tau)_\tau$ at $\alpha$ and $B$ is a standard Brownian motion. We define, for all $0 \leq s \leq t \leq T$,

$$M_s := \exp\left(\int_0^s \psi_{t-r}(x_r)\,dr\right) V_{t-s}(x_s).$$

By Itô's formula and (3.15), we find that

$$dM_s = \exp\left(\int_0^s \psi_{t-r}(x_r)\,dr\right) \psi_{t-s}(x_s) F_{t-s}(x_s)\,ds + dm_s,$$

where $m$ is a martingale. Integrating over $s \in [0,t]$ and taking expectations, we obtain

$$V_t(x) = \mathbf{E}_x\left[\int_0^t \exp\left(\int_0^s \psi_{t-r}(x_r)\,dr\right) \psi_{t-s}(x_s) F_{t-s}(x_s)\,ds\right].$$

Using the hypothesis $\psi \leq 0$, we find that

$$|V_t(x)| \leq -\|F\|_T \mathbf{E}_x\left[\int_0^t \exp\left(\int_0^s \psi_{t-r}(x_r)\,dr\right) \psi_{t-s}(x_s)\,ds\right].$$

The $ds$ integral inside the expectation can be evaluated explicitly and equals

$$\exp\left(\int_0^t \psi_{t-r}(x_r)\,dr\right) - 1 \geq -1.$$

Therefore, $|V_t(x)| \leq \|F\|_T$ and (3.16) is proved in the case of Neumann boundary conditions. The case of Dirichlet boundary conditions follows similarly by killing $(x_\tau)_{\tau \geq 0}$ if it hits $b$ or $c$ before time $t$. $\square$

PROOF OF LEMMA 3.1. We recall that the solutions of (3.1) and (3.2) are constructed in [18], respectively, [11], as monotone nondecreasing limits for $\epsilon \downarrow 0$ and $\lambda \downarrow 0$ of solutions $z = z^{\epsilon,\lambda,\delta}$ of the s.p.d.e.

$$
(3.17) \quad \begin{cases} \dfrac{\partial z}{\partial t} = \dfrac{1}{2}\dfrac{\partial^2 z}{\partial x^2} + f_{\epsilon,\lambda,\delta}(z) + \dfrac{\partial^2 W}{\partial t\,\partial x}, \\ z_t(b) = z_t(c) = a, \qquad t \geq 0, \\ z_0(x) = \overline{v}(x), \qquad x \in [b,c], \end{cases}
$$

where $f_{\epsilon,\lambda,\delta} := f_1 + f_2$ and using the notation (1.4),

$$f_1(r) := \frac{\arctan([r \wedge 0]^2)}{\epsilon},$$

$$f_2(r) := \frac{c_\delta}{\lambda + [r \vee 0]^3}, \qquad r \in \mathbb{R},$$

and $\epsilon, \lambda > 0$. Notice that [11] and [18] use $f_1(r) = r^-/\epsilon$ instead of the definition above: our choice does not change the limit of $z^{\epsilon,\lambda,\delta}$ as $\epsilon \downarrow 0$ and $\lambda \downarrow 0$, but it makes $f_1(\cdot)$ differentiable at 0.

Observe that for fixed $\epsilon < \epsilon'$, $\lambda < \lambda'$ and $\delta > \delta' \geq 3$, Theorem I.3.1 of [12] implies that $z^{\epsilon,\lambda,\delta} \geq z^{\epsilon',\lambda',\delta'}$ and, therefore, $c \leq c'$ implies $u^{(c)} \leq u^{(c')}$.

Therefore, it is enough to prove that there exist finite random variables $\gamma_1$ and $\gamma_2$, independent of $\epsilon, \lambda > 0$, such that

$$(3.18) \quad |z_t^{\epsilon,\lambda,\delta}(x) - z_t^{\epsilon,\lambda,\delta}(y)| \leq \gamma_1 |x - y|^\beta, \qquad x, y \in [0,1], \qquad T \geq t \geq 0,$$

and

$$(3.19) \quad z_t^{\epsilon,\lambda,\delta}(x) - z_s^{\epsilon,\lambda,\delta}(x) \geq -\gamma_2 (t-s)^{\beta/2}, \qquad T \geq t \geq s \geq 0, \qquad x \in [0,1].$$

Notice that $f_{\epsilon,\lambda,\delta}$ is nonnegative and bounded with (bounded) Lipschitz-continuous derivative $f'_{\epsilon,\lambda,\delta}$ (the bounds depend on $\epsilon$, $\lambda$, $\delta$), and $f'_{\epsilon,\lambda,\delta} \leq 0$ over $\mathbb{R}$. Moreover, since either (I) or (II) in (3.3) is satisfied, for $\delta = 3$ and $r \geq 0$ or for $\delta > 3$ and all $r > 0$, we have

$$(3.20) \quad \sup_{\epsilon,\lambda} f_{\epsilon,\lambda,\delta}(r) < \infty.$$

PROOF OF (3.18). For $\eta \in (0,T)$, set $O_{T,\eta} := [\eta, T] \times [b,c]$, and for $\beta \in (0,1)$, denote by $C^{\beta/2,\beta}(O_{T,\eta})$ the set of continuous $N : O_{T,\eta} \mapsto \mathbb{R}$ such that

$$[N]_{\beta/2,\beta} := \sup_{\eta < s < t < T} \sup_{b < x < y < c} \frac{|N_t(x) - N_s(y)|}{|t-s|^{\beta/2} + |x-y|^\beta} < \infty.$$

Moreover, let $C_0^{\beta/2,\beta}(O_{T,\eta})$ be the set of all $N \in C^{\beta/2,\beta}(O_{T,\eta})$ such that $N_t(b) = N_t(c) = 0$ for all $t \in [\eta, T]$. When $\eta = 0$, we write $O_T$ instead of $O_{T,0}$.

It is easy to check that $z_t^{\epsilon,\lambda,\delta}(x) = a + w_t^{(S)}(x) + S_t^{(\bar{v}-a)}(x)$, where for any function $N \in C_0^{\beta/2,\beta}(O_T)$, $w = w^{(N)}$ is the unique solution of the partial differential equation (p.d.e.)

$$\begin{cases} \dfrac{\partial w_t(x)}{\partial t} = \dfrac{1}{2}\dfrac{\partial^2 w_t(x)}{\partial x^2} + f_{\epsilon,\lambda,\delta}(a + w_t(x) + N_t(x)), \\ w_0(x) = 0, & x \in [b,c], \\ w_t(b) = w_t(c) = 0, & t \geq 0, \end{cases}$$

and $S = S^{(\bar{v}-a)}$ is defined in (3.6), with $\bar{v}$ replaced by $\bar{v} - a$. Clearly, $w^{(N)} = h^{\epsilon,\lambda,\delta} + k^{(N)}$, where for all $N \in C_0^{\beta/2,\beta}(O_T)$, $h = h^{\epsilon,\lambda,\delta}$ and $k = k^{(N)}$ are the unique solutions of

(3.21) $$\begin{cases} \dfrac{\partial h_t(x)}{\partial t} = \dfrac{1}{2}\dfrac{\partial^2 h_t(x)}{\partial x^2} + f_{\epsilon,\lambda,\delta}(a), \\ h_0(x) = 0, & x \in [b,c], \\ h_t(b) = h_t(c) = 0, & t \geq 0, \end{cases}$$

and

(3.22) $$\begin{cases} \dfrac{\partial k^{(N)}}{\partial t} = \dfrac{1}{2}\dfrac{\partial^2 k^{(N)}}{\partial x^2} + f_{\epsilon,\lambda,\delta}(a + k^{(N)} + h + N) - f_{\epsilon,\lambda,\delta}(a), \\ k_0^{(N)}(x) = 0, & x \in [b,c], \\ k_t^{(N)}(b) = k_t^{(N)}(c) = 0, & t \geq 0. \end{cases}$$

Express $h_t(x)$ as the convolution of the Green function $g$ and the constant $f_{\epsilon,\lambda,\delta}(a)$, and use (3.20) and the integrability of the partial derivative of $g$ with respect to $x$ to see that

(3.23) $$\sup_{\epsilon,\lambda>0} \|\partial_x h^{\epsilon,\lambda,\delta}\|_T = \kappa(a,\delta,T) < \infty,$$

where $\|\cdot\|_T$ denotes the sup-norm over $O_T$.

Fix $N, M \in C_0^{\beta/2,\beta}(O_T)$ and set $V := k^{(N)} - k^{(M)}$. Then by the mean value theorem, we find that $V$ satisfies (3.15) with $F := N - M$ and

$$\psi_t(x) = f'_{\epsilon,\lambda,\delta}(r_t(x)) \leq 0,$$

where $r_t(x)$ is some number between $a + k_t^{(N)}(x) + h_t(x) + N_t(x)$ and $a + k_t^{(M)}(x) + h_t(x) + M_t(x)$. Moreover, $V$ satisfies homogeneous Dirichlet boundary conditions. By Lemma 3.6, we obtain

(3.24) $$\|k^{(N)} - k^{(M)}\|_T \leq \|N - M\|_T, \qquad N, M \in C_0^{\beta/2,\beta}(O_T),$$

where $\|\cdot\|_T$ denotes the sup-norm over $O_T$. We notice that the same estimate can also be proven with the arguments of [11], (B), page 83.

We now claim that for each $\beta \in (0,1)$,

$$\text{(3.25)} \qquad \sup_{\epsilon,\lambda>0} \sup_{0<t<T} \sup_{b<x<y<c} \frac{|k_t^{(S)}(x) - k_t^{(S)}(y)|}{|x-y|^\beta} < \infty.$$

To establish this, notice first, by [6], Proposition 7.3.2, that $k^{(N)} \in C^{1,2}(O_{T,\eta})$ and that the inhomogeneous term in (3.22) vanishes at $x = b$ and $x = c$. Since $\frac{\partial k^{(N)}}{\partial t}(x) = 0$ for $x \in \{b,c\}$, we see by continuity that

$$\text{(3.26)} \qquad \frac{\partial^2 k_t^{(N)}}{\partial x^2}(b) = \frac{\partial^2 k_t^{(N)}}{\partial x^2}(c) = 0, \qquad t \in (0,T].$$

Recall that $S = S^{(\bar{v}-a)}$ is the stochastic convolution defined by (3.6) above, with $\bar{v}$ replaced by $\bar{v} - a$. For $\rho > 0$, set

$$S_t^\rho(x) := \int_b^c g_{\rho^2}(x,y) S_t(y)\, dy, \qquad x \in [b,c], \qquad t \geq 0.$$

By (3.7), $S$ belongs to $C_0^{\beta/2,\beta}(O_T)$. Therefore, $S^\rho$ belongs to $C_0^{\beta/2,\beta}(O_T)$ and admits a partial derivative in $x$, $\partial_x S^\rho \in C^{\beta/2,\beta}(O_T)$. Moreover, a direct calculation shows that there exists a constant $C_\beta < \infty$ such that a.s. for all $\rho > 0$,

$$\text{(3.27)} \qquad \|S - S^\rho\|_T \leq C_\beta \rho^\beta \gamma_S, \qquad \|\partial_x S^\rho\|_T \leq \frac{C_\beta}{\rho^{1-\beta}} \gamma_S,$$

where $\gamma_S$ is the random variable in (3.7) above. In particular, by (3.24) and (3.27),

$$\text{(3.28)} \qquad \|k^{(S)} - k^{(S^\rho)}\|_T \leq C_\beta \gamma_S \rho^\beta.$$

Let $\tilde{w}$ be the solution of the p.d.e.

$$\begin{cases} \dfrac{\partial \tilde{w}}{\partial t} = \dfrac{1}{2} \dfrac{\partial^2 \tilde{w}}{\partial x^2} + f'_{\epsilon,\lambda,\delta}(a + k^{(S^\rho)} + h^{\epsilon,\lambda,\delta} + S^\rho) \cdot (\tilde{w} + \partial_x h^{\epsilon,\lambda,\delta} + \partial_x S^\rho), \\ \tilde{w}_0(x) = 0, \qquad x \in [b,c], \\ \dfrac{\partial \tilde{w}_t}{\partial x}(b) = \dfrac{\partial \tilde{w}_t}{\partial x}(c) = 0, \qquad t \geq 0. \end{cases}$$

Choosing $N = S^\rho$ and formally differentiating (3.22) with respect to $x$, we see that, in fact, $\tilde{w} = \partial_x k^{(S^\rho)}$ [note that the boundary conditions for $\tilde{w}$ are compatible with (3.26)]. Moreover, setting $V := \tilde{w}$, then $V$ satisfies (3.15) with

$$\psi := f'_{\epsilon,\lambda,\delta}(a + h^{\epsilon,\lambda,\delta} + k + S^\rho) \leq 0, \qquad F := \partial_x h^{\epsilon,\lambda,\delta} + \partial_x S^\rho,$$

and with homogeneous Neumann boundary conditions. Therefore, by Lemma 3.6,

$$\text{(3.29)} \qquad \|\partial_x k^{(S^\rho)}\|_T \leq \|\partial_x h^{\epsilon,\lambda,\delta}\|_T + \|\partial_x S^\rho\|_T.$$

Therefore, by (3.23), (3.27), (3.28) and (3.29), there exists a finite random variable $\gamma_k$, not depending on $\epsilon$ or $\lambda$, such that

$$(3.30) \qquad \|k^{(S)} - k^{(S^\rho)}\|_T \leq \gamma_k \rho^\beta, \qquad \|\partial_x k^{(S^\rho)}\|_T \leq \frac{\gamma_k}{\rho^{1-\beta}}.$$

It follows that

$$(3.31) \ |k_t^{(S)}(x) - k_t^{(S)}(y)| \leq 3\gamma_k |x-y|^\beta, \qquad x, y \in [b, c], \qquad t \in [0, T].$$

Indeed, for $x, y \in [b, c]$, setting $\rho := |x - y|$, by (3.30),

$$\begin{aligned}
|k_t^{(S)}(x) &- k_t^{(S)}(y)| \\
&\leq |k_t^{(S)}(x) - k_t^{(S^\rho)}(x)| \\
&\quad + |k_t^{(S^\rho)}(x) - k_t^{(S^\rho)}(y)| + |k_t^{(S^\rho)}(y) - k_t^{(S)}(y)| \\
&\leq 2\gamma_k \rho^\beta + \frac{\gamma_k}{\rho^{1-\beta}}|x-y| = 3\gamma_k|x-y|^\beta.
\end{aligned}$$

Since (3.31) is uniform in $\epsilon, \lambda$, we obtain (3.25).

By (3.23) and (3.25), we obtain (3.18) with $\gamma_1 := \kappa(a, \delta, T) + 3\gamma_k + \gamma_S$. □

PROOF OF (3.19). The mild formulation of (3.17) yields

$$\begin{aligned}
z_t(x) = \int_b^c g_{t-s}(x, y) z_s(y) \, dy &+ \int_s^t \int_b^c g_{t-r}(x, y) f_{\epsilon, \lambda, \delta}(z_r(y)) \, dy \, dr \\
&+ \int_s^t \int_b^c g_{t-r}(x, y) W(dr, dy).
\end{aligned}$$

Since $f_{\epsilon, \lambda, \delta} \geq 0$, by (3.7),

$$z_t(x) - z_s(x) \geq -\int_b^c g_{t-s}(x, y)|z_s(y) - z_s(x)| \, dy - \gamma_S(t-s)^{\beta/2}$$

for all $T \geq t \geq s \geq 0$. By (3.18) and a standard Gaussian estimate for $g$,

$$\begin{aligned}
\int_b^c g_{t-s}(x, y)|z_s(y) - z_s(x)| \, dy &\leq \int_\mathbb{R} \frac{\gamma_{z,1}|y|^\beta}{\sqrt{2\pi(t-s)}} e^{-y^2/(2(t-s))} \, dy \\
&\leq \gamma_{z,1}(t-s)^{\beta/2}.
\end{aligned}$$

Therefore, we obtain (3.19) with $\gamma_2 := \gamma_{z,1} + \gamma_S$. The proof of Lemma 3.1 is complete. □

REMARK 3.7. In the case where $\delta > 6$, one can easily obtain actual Hölder continuity in time of the solution $v$ of (3.1), rather than the lower bound obtained in (3.5). Consider for simplicity the case $[b, c] = [0, 1]$. More precisely, using the mild formulation of (3.1), it suffices to consider the process

$$v_t(x) = S_t^{(\bar{v})}(x) + \int_0^t \int_0^1 g_{t-s}(x,y) \frac{c_\delta}{(v_s(y))^3} \, ds \, dy,$$

where $S^{(\bar{v})}$ is defined in (3.6) with $b=0$ and $c=1$, and $\bar{v}$ is a Bessel bridge of dimension $\delta$ independent of $W$. The first term is Hölder-continuous in $t$ with exponent $1/4$, by (3.7), so we check this property for the second term.

Fix $\epsilon > 0$ and split the $dy$ integral into three integrals, over $[0, \epsilon]$, $[1-\epsilon, 1]$ and $[\epsilon, 1-\epsilon]$, yielding, respectively, three terms $v_t^{(1)}(x)$, $v_t^{(2)}(x)$ and $v_t^{(3)}(x)$. For $x \in [2\epsilon, 1-2\epsilon]$, the first two terms are $C^\infty$. Notice that for such $x$ and $0 < t_1 < t_2 < T$, by the Cauchy–Schwarz inequality,

$$(v_{t_1}^{(3)}(x) - v_{t_2}^{(3)}(x))^2 \leq \left[ \int_0^T \int_\epsilon^{1-\epsilon} (g_{t_1-s}(x,y) \mathbf{1}_{\{s \leq t_1\}} - g_{t_2-s}(x,y) \mathbf{1}_{\{s \leq t_2\}})^2 \, ds \, dy \right]$$
$$\times \int_0^T \int_\epsilon^{1-\epsilon} \frac{c_\delta^2}{(v_s^{(3)}(y))^6} \, ds \, dy.$$

It is well known ([1], Lemma B.1) that the first factor is bounded by $C(t_2 - t_1)^{1/2}$, so it suffices to check that the second factor is finite a.s. Using the explicit form of the marginal densities of the Bessel bridge [14], one checks that this is indeed the case.

**4. Proof of Theorem 2.1.** The proof is based on Theorem 3.5 and on a comparison technique.

For $\delta \geq 3$, we define $u^{(\delta)}$ as follows: $u^{(3)}$ is the solution of (2.2) and, for $\delta > 3$, $u^{(\delta)}$ is the solution of (2.1). Notice that the result concerning $u^{(3)}$ is already established by Theorem 3.5, since (I) in (3.3) and (3.12) are satisfied by $k = 5$ and $\delta = 3$.

Let $\delta > 3$. By the monotonicity in $\delta$ (see the beginning of the proof of Lemma 3.1), almost surely, $u^{(\delta)} \geq u^{(3)}$. Let $\epsilon \in (0, 1/2)$, $T > 0$ and $\beta \in (1/4, 1/2)$. Consider the intervals $I_1 = [0, \epsilon]$ and $I_2 = [1-\epsilon, 1]$, and the random variable

$$\eta := \inf_{t \in [\epsilon, T]} \min_{i=1,2} \sup_{x \in I_i} u_t^{(3)}(x).$$

By the result just established for $u^{(3)}$, there is no $t \in [\epsilon, 1]$ such that $u_t^{(3)}(\cdot)$ vanishes identically on $I_1$ or $I_2$, and therefore $\eta > 0$ a.s.

For $n, j \in \mathbb{N}$, set $t_{n,j} = j 2^{-8n}$ and let $j_n = \inf\{j \geq 0 : t_{n,j} > \epsilon\}$. Let $X_{n,j,i}$ be the leftmost (but in fact unique) point in $I_i$ such that

$$u_{t_{n,j}}^{(3)}(X_{n,j,i}) = \sup_{x \in I_i} u_{t_{n,j}}^{(3)}(x).$$

Then $X_{n,j,i}$ is $\mathcal{F}_{t_{n,j}}$-measurable.

Let $\gamma_v$ be the random variable that appears in (3.5) (for $\delta = 3$, $a = 0$). For $n \in \mathbb{N}$, let

$$F_n = \{\eta > 2^{-n+1}, \gamma_v < 2^{n(4\beta-1)}\}.$$

Because $\eta > 0$ a.s. and $\gamma_v < \infty$ a.s.,

$$P\left(\bigcup_{n \in \mathbb{N}} F_n\right) = 1.$$

We claim that for all $n \in \mathbb{N}$, $\omega \in F_n$, $j \in [j_n, T2^{4n}]$, $t \in [t_{n,j}, t_{n,j+1}]$ and $x = X_{n,j,1}$ or $x = X_{n,j,2}$,

(4.1) $$u_t^{(\delta)}(x) \geq 2^{-n}.$$

Indeed, $u^{(\delta)} \geq u^{(3)}$ a.s. Moreover, by the definition of $F_n$, for $\omega \in F_n$, we have

$$u_{t_{n,j}}^{(3)}(X_{n,j,i}) \geq \eta > 2^{-n+1}, \qquad i = 1, 2.$$

Also, for $t \in [t_{n,j}, t_{n,j+1}]$, by (3.5),

$$u_t^{(3)}(x) - u_{t_{n,j}}^{(3)}(x) \geq -\gamma_v (t - t_{n,j})^{\beta/2} \geq -\gamma_v 2^{-4\beta n} \geq -2^{-n}.$$

Finally, for $\omega \in F_n$ and $t \in [t_{n,j}, t_{n,j+1}]$,

(4.2) $$u_t^{(\delta)}(X_{n,j,i}) \geq u_t^{(3)}(X_{n,j,i}) \geq u_{t_{n,j}}^{(3)}(X_{n,j,i}) - 2^{-n} > 2^{-n}, \qquad i = 1, 2.$$

Now let $\tilde{u}$ be the solution of (3.1) in the domain $[t_{n,j}, t_{n,j+1}] \times [b, c]$, where $b = X_{n,j,1}$ and $c = X_{n,j,2}$, with (random) initial condition

$$\tilde{u}_{t_{n,j}}(\cdot) = \begin{cases} 2^{-n}, & \text{if } \min(u_{t_{n,j}}^{(\delta)}(b), u_{t_{n,j}}^{(\delta)}(c)) \leq 2^{-n}, \\ \min(u_{t_{n,j}}^{(\delta)}(\cdot), 2^{-n}), & \text{otherwise}, \end{cases}$$

and boundary conditions $2^{-n}$. Notice that Theorem 3.5 applies to $\tilde{u}$: since the initial condition is $\mathcal{F}_{t_{n,j}}$-measurable, we can condition on this $\sigma$-field.

We claim that the following holds:

(4.3) $$u_t^{(\delta)}(x) \geq \tilde{u}_t(x), \qquad (t, x) \in [t_{n,j}, t_{n,j+1}] \times [b, c], \qquad \omega \in F_n.$$

Since $\tilde{u}$ has the desired property by Theorem 3.5, it would follow that $u^{(\delta)}$ does too. Thus, (4.3) would finish the proof of Theorem 2.1.

To establish (4.3), we consider again the process $z = z^{\epsilon, \lambda, \delta}$, which solves the s.p.d.e. (3.17) with $a$ and $b$ replaced by 0, $c$ replaced by 1 and $\bar{v}$ replaced by $\bar{u}$. Recall that $u^{(\delta)}$ is the monotone limit of $z^{\epsilon, \lambda, \delta}$ as $\epsilon \searrow 0$ and then $\lambda \searrow 0$. In particular,

$$u_t^{(\delta)}(x) = \sup_{\epsilon > 0, \lambda > 0} z_t^{\epsilon, \lambda, \delta}(x).$$

For $\omega \in F_n$ and $t \in [t_{n,j}, t_{n,j+1}]$, by (4.2), $u_t^{(\delta)}(b) > 2^{-n}$ and $u_t^{(\delta)}(c) > 2^{-n}$. By Dini's theorem, the convergence of $z^{\epsilon,\lambda,\delta}$ to $u^{(\delta)}$ is uniform on $O_T$, so that we can find $\Theta(\omega)$ such that for all $\epsilon \leq \Theta(\omega)$ and $\lambda \leq \Theta(\omega)$, we have $z_t^{\epsilon,\lambda,\delta}(b) > 2^{-n}$ and $z_t^{\epsilon,\lambda,\delta}(c) > 2^{-n}$.

For all such $\epsilon, \lambda$, let $\tilde{z} = \tilde{z}^{\epsilon,\lambda,\delta}$ be the solution of (3.17) in the domain $[t_{n,j}, t_{n,j+1}] \times [b, c]$, where $b = X_{n,j,1}$ and $c = X_{n,j,2}$, with (random) initial condition

$$\tilde{z}_{t_{n,j}}^{\epsilon,\lambda,\delta}(\cdot) = \begin{cases} 2^{-n}, & \text{if } \min(z_{t_{n,j}}^{\epsilon,\lambda,\delta}(b), z_{t_{n,j}}^{\epsilon,\lambda,\delta}(c)) \leq 2^{-n}, \\ \min(z_{t_{n,j}}^{\epsilon,\lambda,\delta}(\cdot), 2^{-n}), & \text{otherwise}, \end{cases}$$

and boundary conditions $a := 2^{-n}$. Setting $V_t(x) := z_{t+t_{n,j}}^{\epsilon,\lambda,\delta}(x) - \tilde{z}_{t+t_{n,j}}^{\epsilon,\lambda,\delta}(x)$, by the mean value theorem, we have

$$\frac{\partial V}{\partial t} = \frac{1}{2}\frac{\partial^2 V}{\partial x^2} + f_{\epsilon,\lambda,\delta}(z) - f_{\epsilon,\lambda,\delta}(\tilde{z}) = \frac{\partial^2 V}{\partial x^2} + \psi \cdot V,$$

where $\psi : O_T \mapsto \mathbb{R}$ is bounded and $\psi \leq 0$. Moreover, on $F_n$,

$$V_0(x) \geq 0, \qquad V_t(b) \geq 0, \qquad V_t(c) \geq 0.$$

Since $f_{\epsilon,\lambda,\delta}(z) - f_{\epsilon,\lambda,\delta}(\tilde{z})$ is in $C_0^{\beta/2,\beta}(O_T)$, it follows that $V$ is in $C^{1,2}(O_T)$ and we can apply the maximum principle (see, e.g., [13], Chapter 3, Theorem 7 and the remark on page 174) to obtain $V \geq 0$. In particular, on $F_n$, the following holds:

$$z^{\epsilon,\lambda,\delta} \geq \tilde{z}^{\epsilon,\lambda,\delta} \qquad \text{on } [t_{n,j}, t_{n,j+1}] \times [b,c] \qquad \text{for all } \epsilon, \lambda \leq \Theta(\omega).$$

Taking $\epsilon, \lambda \to 0$, we get (4.3), which finishes the proof of Theorem 2.1.

**5. The Gaussian random string.** In this section, we prove Theorem 2.4. We prove first the "positive" assertions not contained in Theorem 2.3, which we summarize in the following lemma.

LEMMA 5.1.

- If $d = 3$, then with positive probability, there exist $t > 0$ and $x_1 < x_2 < x_3$ such that $U_t(x_i) = 0$, $i = 1, 2, 3$.
- If $d = 2$, then for all $k \in \mathbb{N}$, with positive probability, there exist $t > 0$ and $x_1 < \cdots < x_k$ such that $U_t(x_i) = 0$, $i = 1, \ldots, k$.

We recall the following results, proved, respectively, in Proposition 1 and Corollary 3 of [9].

PROPOSITION 5.2. *The components $U^i$ of the stationary pinned string are mutually independent centered Gaussian random fields with covariance function determined by*

$$\mathbb{E}[(U_t^i(x) - U_s^i(y))^2] =: c(t, x; s, y),$$

*where $c$ is such that there is a constant $c_1 > 0$ such that for all $x, y \in \mathbb{R}$ and $0 \le s \le t$,*

$$c_1(|x - y| + |t - s|^{1/2}) \le c(t, x; s, y) \le 2(|x - y| + |t - s|^{1/2}).$$

PROPOSITION 5.3. *For any compact set $A \subset (0, \infty) \times \mathbb{R}$, the laws of the random fields $(U_t(x) : (t, x) \in A)$ and $(U_t(x) + z : (t, x) \in A)$ are mutually absolutely continuous.*

The following result is stated in Corollary 5 of [9] (see also [10]).

PROPOSITION 5.4. *For any compact sets $A^+ \subset (0, \infty) \times (0, \infty)$ and $A^- \subset (0, \infty) \times (-\infty, 0)$, the law of*

$$((U_t(x) : (t, x) \in A^+), (U_t(x) : (t, x) \in A^-))$$

*and the law of*

$$((U_t(x) : (t, x) \in A^+), (\tilde{U}_t(x) : (t, x) \in A^-))$$

*are mutually absolutely continuous, where $U$ and $\tilde{U}$ are independent copies of the stationary pinned string.*

PROOF OF LEMMA 5.1. Let $d = 3$ and $k = 3$ or let $d = 2$ and $k \in \mathbb{N}$, and for $t \in [1, 2]$ and $x_i \in [2i, 2i + 1]$, $i = 1, \ldots, k$, set

$$Z(t, x_1, \ldots, x_k) := (U_t(x_1), \ldots, U_t(x_k)) \in \mathbb{R}^{kd},$$

$$X(t, x_1, \ldots, x_k) := (\tilde{U}_t^{(1)}(x_1), \ldots, \tilde{U}_t^{(k)}(x_k)) \in \mathbb{R}^{kd},$$

where $(\tilde{U}^{(i)})_i$ are i.i.d. copies of the stationary pinned string in $\mathbb{R}^d$. The lemma will follow if we prove that $0 \in \mathbb{R}^{kd}$ belongs to the range of $Z$ with positive probability. By Proposition 5.4, the laws of $Z$ and $X$ are mutually absolutely continuous. Therefore, it suffices to prove that $0 \in \mathbb{R}^{kd}$ belongs to the range of $X$ with positive probability.

We use results on existence of occupation densities for Gaussian processes proved in Sections 6 and 22 of [2]. Let $T := [1, 2] \times \prod_{i=1}^{k} [2i, 2i + 1]$ and, for $\tau = (t, x_1, \ldots, x_k) \in T$, set $X_\tau := X(t, x_1, \ldots, x_k) \in \mathbb{R}^{kd}$. Then $X : T \mapsto \mathbb{R}^{kd}$ is a centered continuous Gaussian process such that the determinant $\Delta(\sigma, \tau)$ of the covariance matrix of $X_\sigma - X_\tau$ is positive for almost all $(\sigma, \tau) \in T \times T$. Following [2], we say that $X$ is (LT) (short for local time) if there exists a (random) measurable kernel $(\alpha(z, A) : z \in \mathbb{R}^d, A \subseteq T$ Borel), termed *occupation kernel*, such that a.s., for all bounded Borel $f$ over $\mathbb{R}^d$ and $A \subseteq T$ Borel,

$$\int_A f(X_\tau) \, d\tau = \int_{\mathbb{R}^d} f(z) \alpha(z, A) \, dz.$$

We recall [2], Theorem 6.4(ii), that a.s., for all $z \in \mathbb{R}^d$,

(5.1) $\quad \alpha(z, T \setminus M_z) = 0 \quad$ where $M_z := \{\tau \in T : X_\tau = z\}$.

We want to show that $X$ is (LT). As proved in Theorem 22.1 of [2],

(5.2) $\quad \sup_{\sigma \in T} \int_T (\Delta(\sigma, \tau))^{-1/2} \, d\tau < \infty \quad \Longrightarrow \quad X$ is (LT),

so we check that $\Delta(\sigma, \tau)$ has this property. By the independence of the coordinates and Proposition 5.2, for some constant $C > 0$,

$$\Delta(\sigma, \tau) = \prod_{i=1}^{k} [c(t, x_i; s, y_i)]^d \geq C \prod_{i=1}^{k} (|x_i - y_i| + |t - s|^{1/2})^d$$

for all $\tau = (t, x_1, \ldots, x_k)$ and $\sigma = (s, y_1, \ldots, y_k) \in T$.

If $d = 3$ and $k = 3$, then

$$\int_T (\Delta(\sigma, \tau))^{-1/2} \, d\tau \leq C \int_0^1 dt \left[ \int_0^1 dx \, (x + t^{1/2})^{-3/2} \right]^3$$

$$\leq C \int_0^1 \frac{1}{t^{3/4}} \, dt < \infty.$$

If $d = 2$, then for any $k \in \mathbb{N}$,

$$\int_T (\Delta(\sigma, \tau))^{-1/2} \, d\tau \leq C \int_0^1 dt \left[ \int_0^1 dx \, (x + t^{1/2})^{-1} \right]^k$$

$$= C \int_0^1 [\log(t^{-1/2} + 1)]^k \, dt < \infty.$$

Therefore, $\alpha(\cdot, \cdot)$ is well defined for such values of $d$ and $k$.

For all bounded Borel $f : \mathbb{R}^d \mapsto \mathbb{R}$,

$$\mathbb{E}\left[\int_T f(X_\tau) \, d\tau\right] = \int_{\mathbb{R}^d} f(z) \mathbb{E}[\alpha(z, T)] \, dz,$$

so for some $z_0 \in \mathbb{R}^d$, with positive probability, $\alpha(z_0, \cdot)$ is not identically zero. By (5.1), $M_{z_0}$ is nonempty with positive probability. By Proposition 5.3, the laws of $M_{z_0} = \{\tau \in T : X_\tau - z_0 = 0\}$ and $M_0 = \{\tau \in T : X_\tau = 0\}$ are mutually absolutely continuous, so that $M_0$ is nonempty with positive probability and the proof is complete. $\square$

We turn now to the "negative" assertions of Theorem 2.4: those not already given in Theorem 2.3 are summarized in the following lemma.

LEMMA 5.5. *If $d \geq 4$, then the probability that there exist $t \geq 0$ and $x_1 < x_2$ such that $U_t(x_i) = 0$, $i = 1, 2$, is 0. If $d = 3$, then the probability that there exist $t \geq 0$ and $x_1 < \cdots < x_4$ such that $U_t(x_i) = 0$, $i = 1, \ldots, 4$, is 0.*

We recall the following scaling lemma, proven in Corollary 1 of [9].

COROLLARY 5.6 ([9]). *The stationary pinned string has the following properties:*

(1) *Translation invariance. For any $t_0 \geq 0$ and $x_0 \in \mathbb{R}$, the field*
$$(U_{t_0+t}(x_0 \pm x) - U_{t_0}(x_0) : x \in \mathbb{R}, t \geq 0)$$
*has the same law as the stationary pinned string.*
(2) *Scaling. For $L > 0$, the field*
$$(L^{-1} U_{L^4 t}(L^2 x) : x \in \mathbb{R}, t \geq 0)$$
*has the same law as the stationary pinned string.*
(3) *Time reversal. For any $T > 0$, the field*
$$(U_{T-t}(x) - U_T(0) : x \in \mathbb{R}, 0 \leq t \leq T)$$
*has the same law as the stationary pinned string over the interval $[0, T]$.*

PROOF OF LEMMA 5.5. All of these proofs follow the "replication" idea that can be found in [9], which originated in the work of Lévy (see, in particular, [9], Section 4, and [5], Theorem 2.2).

*Case* I: $d \geq 4$. By projecting onto the first four coordinates, we see that it is enough to consider the case $d = 4$.

Note at fixed times such as $t = 0$, that since $(U_0(x), U_0(-x) : x \geq 0)$ are independent four-dimensional Brownian motions indexed by $x$, standard properties of Brownian motion imply that, with probability 1, there do not exist points $x_1 \neq x_2$ such that $U_0(x_1) = U_0(x_2) = 0$ (see, e.g., [5], Theorem 1.1).

Let
$$Z_0(t, x, y) = (U_t(x), U_t(y)).$$
Using Corollary 5.6, we see that it is enough to prove that, with probability 1, there are no points

(5.3) $\qquad t \in [1, 2], \qquad x \in [1, 2], \qquad y \in [-2, -1]$

such that $Z_0(t, x, y) = 0$. By Proposition 5.4, it suffices to consider
$$Z(t, x, y) = (U_t(x), \tilde{U}_t(y))$$
instead of $Z_0$.

For $z \in \mathbb{R}^{2d}$, let $Q(z)$ be the event that there do not exist points $t, x$ and $y$ that satisfy (5.3) such that $Z(t, x, y) = z$. Applying Proposition 5.3 to both $U_t(x)$ and $\tilde{U}_t(y)$, we conclude that for $z \in \mathbb{R}^{2d}$,

$$\mathbb{P}(Q(0)) = 0 \quad \Longleftrightarrow \quad \mathbb{P}(Q(z)) = 0.$$

Next, for $A \subset (0, \infty) \times \mathbb{R}^2$, let $Z(A)$ be the range of $Z(t, x, y)$ for $(t, x, y) \in A$ and let $m(\cdot)$ be Lebesgue measure. By Fubini's theorem,

$$\mathbb{E}[m(Z(A))] = \int_{\mathbb{R}^4} \mathbb{P}(Z(t, x, y) = z \text{ for some } (t, x, y) \in A) \, dz.$$

Therefore, it suffices to show that for $A = (1, 2] \times (1, 2] \times (-2, -1]$,

(5.4) $$\mathbb{E}[m(Z(A))] = 0.$$

To prove (5.4), following Lévy, we use scaling to relate $\mathbb{E}[m(Z(A))]$ to $\mathbb{E}[m(Z(A_i))]$, where $A_i$ are certain subsets of $A$. We then show that $\mathbb{E}[m(Z(A_i) \cap Z(A_j))] = 0$ for $i \neq j$. An independence argument will imply that $\mathbb{E}[m(Z(A_i))] = 0$ for each $i$ and so $\mathbb{E}[m(Z(A))] = 0$, and we will be finished.

Subdivide the cube $A$ into 16 pairwise disjoint subsets as follows. Subdivide each space interval $[1, 2]$ and $[-2, -1]$ into two disjoint subintervals of equal length and subdivide the time interval $[1, 2]$ into four disjoint subintervals of equal length. All these subintervals are taken open on the left and closed on the right. By taking cartesian products, form 16 disjoint sets $A_i$, $i = 1, \ldots, 16$, whose union is $A$.

Now we use Corollary 5.6 to scale time and space: we find that

$$(Z(4t, 2x, 2y), (t, x, y) \in \mathbb{R}^3) \stackrel{\mathcal{D}}{=} (2^{1/2} Z(t, x, y), (t, x, y) \in \mathbb{R}^3),$$

where the equality is in distribution. Since $Z(t, x, y)$ is a vector with two coordinates, each of which lies in $\mathbb{R}^4$, the range of $Z(t, x, y)$ lies in $\mathbb{R}^8$. Therefore,

$$\mathbb{E}[m(Z(A))] = (2^{1/2})^8 \mathbb{E}[m(Z(A_i))] = 16 \mathbb{E}[m(Z(A_i))].$$

A standard inclusion–exclusion argument implies

$$\mathbb{E}[m(Z(A))] \leq \sum_{i=1}^{16} \mathbb{E}[m(Z(A_i))] - \sum_{1 \leq i < j \leq 16} \mathbb{E}[m(Z(A_i) \cap Z(A_j))].$$

Therefore, for each pair $i < j$,

$$\mathbb{E}[m(Z(A_i) \cap Z(A_j))] = 0.$$

Next, relabelling if necessary, choose $A_1, A_2$ such that the points in $A_1, A_2$ have the same $x, y$ coordinates, but the $t$ coordinates lie in adjacent time intervals. Let $t_0$ be the common boundary point of these two time intervals.

Let $\mathcal{H}$ be the $\sigma$-field generated by the values of $U_t(\cdot), \tilde{U}_t(\cdot)$ for $t = t_0$. Note that $\mathcal{H}$ is also generated by the values of $Z(t, x, y)$ for $t = t_0$. By the Markov property in time of $Z(\cdot, \cdot, \cdot)$, the random variables $Z(A_1)$ and $Z(A_2)$ are conditionally independent given $\mathcal{H}$, and by the time-reversal property of $U_t(\cdot)$ and $\tilde{U}_t(\cdot)$ given in Corollary 5.6, their conditional distributions given $\mathcal{H}$ coincide. Therefore, using versions of conditional expectations that are jointly measurable in $(z, \omega)$ ([15], Lemma 3), we obtain

$$\begin{aligned}
0 &= \mathbb{E}[m(Z(A_1) \cap Z(A_2))] \\
&= \int_{\mathbb{R}^4} \mathbb{E}[\mathbf{1}_{\{z \in Z(A_1)\}} \mathbf{1}_{\{z \in Z(A_2)\}}] \, dz \\
&= \mathbb{E}\left( \int_{\mathbb{R}^4} \mathbb{E}[\mathbf{1}_{\{z \in Z(A_1)\}} \mathbf{1}_{\{z \in Z(A_2)\}} | \mathcal{H}] \, dz \right) \\
&= \mathbb{E}\left( \int_{\mathbb{R}^4} \mathbb{E}[\mathbf{1}_{z \in Z(A_1)\}} | \mathcal{H}] \mathbb{E}[\mathbf{1}_{\{z \in Z(A_2)\}} | \mathcal{H}] \, dz \right) \\
&= \mathbb{E}\left( \int_{\mathbb{R}^4} \mathbb{E}[\mathbf{1}_{\{z \in Z(A_1)\}} | \mathcal{H}]^2 \, dz \right).
\end{aligned}$$

This implies that $\mathbb{E}[\mathbf{1}_{\{z \in Z(A_1)\}} | \mathcal{H}] = 0$ for almost every $z$, a.s. Therefore,

$$\mathbb{E}[m(Z(A))] = 16 \mathbb{E}[m(Z(A_1))] = 16 \mathbb{E}\left[ \int_{\mathbb{R}^4} \mathbb{E}[\mathbf{1}_{\{z \in Z(A_1)\}} | \mathcal{H}] \, dz \right] = 0$$

and, hence, $m(Z(A)) = 0$ a.s. This proves (5.4) and completes the proof of Case I ($d \geq 4$).

*Case* II: $d = 3$. Since this proof is similar to the previous case, we only outline the main points. We must show that with probability 1, there do not exist $t \geq 0$ and $x_1 < \cdots < x_4$ with $U_t(x_i) = 0$ for $i = 1, \ldots, 4$. As in the previous proof, we assume that $t, x_1, \ldots, x_4$ lie in a bounded set $A$, namely

$$A := \{(t, x_1, \ldots, x_4) : t \in (1, 2], x_i \in (a_i, b_i], i = 1, \ldots, 4\},$$

where $a_1 < b_1 < a_2 < \cdots < b_4$. Let

$$Z_0(t, x_1, \ldots, x_4) = (U_t(x_1), \ldots, U_t(x_4)).$$

We must show that with probability 1, there does not exist $(t, x_1, \ldots, x_4) \in A$ with $Z_0(t, x_1, \ldots, x_4) = 0$.

Let $U_t^{(i)}(x) : i = 1, \ldots, 4$ be independent copies of $U$ and let

$$Z_t(x_1, \ldots, x_4) := (U_t^{(1)}(x_1), \ldots, U_t^{(4)}(x_4)).$$

By Proposition 5.4, it is enough to show that with probability 1, there exists no $(t, x_1, \ldots, x_4) \in A$ with $Z_t(x_1, \ldots, x_4) = 0$.

Once again, Proposition 5.3 implies that we need only show

$$\mathbb{E}[m(Z(A))] = 0,$$

where we have used the same notation as in the previous case. Now we divide each interval $(a_i, b_i]$ into two equally long subintervals and divide $(1, 2]$ into four subintervals of equal length. The products of these intervals give us 64 "rectangles" $A_i$. Once again, scaling implies that for each value of $i$,

$$\mathbb{E}[m(Z(A))] = (2^{1/2})^{12}\mathbb{E}[m(Z(A_i))]$$
$$= 64\mathbb{E}[m(Z(A_i))].$$

Since 64 is also the number of rectangles, we may argue as before, to conclude that for $i < j$,

$$\mathbb{E}[m(Z(A_i) \cap Z(A_i))] = 0.$$

Then we can use the same conditional independence argument as before to conclude that $\mathbb{E}[m(Z(A_i))] = 0$ and, hence, $\mathbb{E}[m(Z(A))] = 0$.

This completes the proof of Case II ($d = 3$), and Lemma 5.5 is proved. □

All statements in Theorem 2.4 have now been proved.

**6. Proof of Theorem 2.2.** For the proof of Theorem 2.2, we need a different approach. We introduce infinite-dimensional capacities related to the processes $u$ and $U$ of Section 2, and we prove that the former is always greater than or equal to the latter. Since sets of positive capacity are hit with positive probability by the associated Markov process, we use the results of Theorems 2.3 and 2.4 on $U$ and we transfer them to $u$.

PROOF OF THEOREM 2.2. Let $d \in \mathbb{N}$ and denote by $(V_t(x) : t \geq 0, x \in [0, 1])$ the $\mathbb{R}^d$-valued continuous process that is the solution of

(6.1) $$\begin{cases} \dfrac{\partial V}{\partial t} = \dfrac{1}{2}\dfrac{\partial^2 V}{\partial x^2} + \dfrac{\partial^2 W_d}{\partial t \, \partial x}, \\ V_t(0) = V_t(1) = 0, \quad t \geq 0, \\ V_0 = \overline{V} \in (C^+)^d. \end{cases}$$

Let $A := [\epsilon, T] \times [\epsilon, 1 - \epsilon]$. We claim that the laws of $(V_t(x) : (t, x) \in A)$ and $(U_t(x) : (t, x) \in A)$ are mutually absolutely continuous, where $U$ is the stationary pinned string (2.4). Indeed, let $\psi : [0, \infty] \times \mathbb{R} \mapsto [0, 1]$ be a $C^\infty$ function with compact support inside $(0, \infty) \times (0, 1)$ such that $\psi \equiv 1$ on $A$. For $x \notin [0, 1]$, set $V_t(x) = 0$ and define

$$Z_t(x) := U_t(x) + \psi_t(x)(V_t(x) - U_t(x)), \qquad (t, x) \in [0, \infty) \times \mathbb{R}.$$

Then $Z \equiv V$ on $A$ and

$$Z_t(x) - U_t(x)$$

$$= \psi_t(x) \int_0^t \int_{\mathbb{R}} [\mathbf{1}_{\{y \in (0,1)\}} g_{t-s}(x,y) - G_{t-s}(x-y)] W(ds, dy)$$

$$+ \psi_t(x) \int_{\mathbb{R}} [\mathbf{1}_{\{y \in (0,1)\}} g_t(x,y) V_0(y) - G_t(x-y) U_0(y)] \, dy.$$

Using the explicit form of $g_{t-s}(x,y)$ (see, e.g., [16]), we notice that the singularity in $g$ cancels with $G$ and, therefore, $Z - U$ is a $C^\infty$ Gaussian process with compact support in $(0, \infty) \times (0, 1)$. It follows that $Z_0(\cdot) \equiv U_0(\cdot)$,

$$\frac{\partial Z_t}{\partial t} = \frac{1}{2} \frac{\partial^2 Z_t}{\partial x^2} + h + \frac{\partial^2 W_d}{\partial t \, \partial x}, \qquad h := \left( \frac{\partial}{\partial t} - \frac{1}{2} \frac{\partial^2}{\partial x^2} \right)(Z - U),$$

and $(h_t(x), \, t \geq 0, \, x \in \mathbb{R})$ is again a continuous Gaussian process, adapted in time to the filtration of $W$, supported on $[0, T] \times [0, 1]$, with variance bounded over $[0, T] \times [0, 1]$. By Lemmas 1 and 2 of [10], the laws of $Z$ and $U$ over $[0, T] \times \mathbb{R}$ are mutually absolutely continuous. Since $Z \equiv V$ on $A$, the claim is proven.

For $\delta > 3$, let $u^{(\delta)}$ denote the solution of (2.1) and let $u^{(3)}$ denote the solution of (2.2). We now recall that the following properties were proved in [18] (see, in particular, Theorems 3 and 5 there):

- For all $\delta \geq 3$, $(u_t^{(\delta)} : \bar{u} \in C^+, t \geq 0)$ is the diffusion associated with the symmetric Dirichlet form with state space $C^+$, defined by

$$W^{1,2}(\pi_\delta) \ni \varphi, \psi \mapsto \mathcal{D}^\delta(\varphi, \psi) := \tfrac{1}{2} \int_K \langle \nabla \varphi, \nabla \psi \rangle \, d\pi_\delta,$$

where $K = \{\bar{u} \in L^2(0,1) : \bar{u} \geq 0\}$, $\nabla$ denotes the Fréchet differential in the Hilbert space $H := L^2(0,1)$ and $\pi_\delta$ is the law of the Bessel bridge $X$.

- For all $d \in \mathbb{N}$, $\overline{V} \mapsto V$ is the diffusion associated with the Dirichlet form $(\Lambda^d, W^{1,2}(\mu_d))$ on $(C^+)^d$, defined by

$$W^{1,2}(\mu_d) \ni F, G \mapsto \Lambda^d(F, G) := \tfrac{1}{2} \int_{H^d} \langle \overline{\nabla} F, \overline{\nabla} G \rangle_{H^d} \, d\mu_d,$$

where $\mu_d$ is the law of a Brownian bridge of dimension $d$ between 0 and 0 over $[0,1]$, and $\overline{\nabla}$ denotes the gradient in $H^d := L^2(0,1; \mathbb{R}^d)$.

- For $d \in \mathbb{N}$, $d \geq 3$, define $\Phi : (C^+)^d \mapsto C^+$ by $\Phi(y)(\tau) := |y(\tau)|$, $\tau \in [0,1]$. Then $\mathcal{D}^d$ is the image of $\Lambda^d$ under the map $\Phi$, that is, $\pi_d$ is the image of $\mu_d$ under $\Phi$ and

(6.2) $\qquad W^{1,2}(\pi_d) = \{\varphi \in L^2(\pi_d) : \varphi \circ \Phi \in W^{1,2}(\mu_d)\},$

(6.3) $\qquad \mathcal{D}^d(\varphi, \psi) = \Lambda^d(\varphi \circ \Phi, \psi \circ \Phi) \qquad \forall \varphi, \psi \in W^{1,2}(\pi_d).$

Formula (6.3) is based on a simple fact, namely that for any $\varphi \in W^{1,2}(\pi_d)$,

$$\varphi \circ \Phi \in W^{1,2}(\mu_d) \quad \text{and} \quad \overline{\nabla}(\varphi \circ \Phi)(y) = \frac{y}{|y|} \nabla \varphi(|y|)$$

for $\mu_d$-a.e. $y$, which implies that

$$\langle \overline{\nabla}(\varphi \circ \Phi), \overline{\nabla}(\psi \circ \Phi) \rangle_{H^d} = \langle \nabla\varphi, \nabla\psi \rangle_H \circ \Phi, \qquad \mu_d\text{-a.s.},$$

since, for all $\tau \in [0,1]$, $y(\tau)/|y(\tau)| \in \mathbb{R}^d$ has Euclidean norm 1. Formula (6.2) is a deeper result, which however we do not need here.

Recall that the $\mathcal{D}$-capacity of a subset of $C^+$ is defined as follows. We set $\mathcal{D}_1 := \mathcal{D} + \langle \cdot, \cdot \rangle_{L^2(\pi_\delta)}$. For $A \subseteq C^+$ open, let

$$\operatorname{Cap}_\mathcal{D}(A) := \inf\{\mathcal{D}_1(\varphi, \varphi) : \varphi \in W^{1,2}(\pi_\delta), \varphi \geq 1, \pi_\delta\text{-a.e. on } A\}.$$

For any $E \subseteq C^+$, let

$$\operatorname{Cap}_\mathcal{D}(E) := \inf\{\operatorname{Cap}_\mathcal{D}(A) : E \subseteq A \subseteq C^+, A \text{ open}\}.$$

The $\Lambda^d$-capacity of subsets of $(C^+)^d$ is defined analogously. Then, by (6.3), for all $E \subseteq C^+$ and $d \in \mathbb{N}$, $d \geq 3$,

(6.4) $$\operatorname{Cap}_{\mathcal{D}^d}(E) \geq \operatorname{Cap}_{\Lambda^d}(\Phi^{-1}(E)).$$

It is now a classical result of potential theory that a set with positive capacity is hit by the associated Markov process with positive probability, and vice versa. For a proof of this statement in infinite-dimensional settings, see Theorems III.2.11(ii) and IV.5.29(i) in [3].

We set

$$E_3 := \{\bar{u} \in C^+ : \exists\, 0 < x_1 < x_2 < x_3 < 1, \bar{u}(x_i) = 0, i = 1, 2, 3\}.$$

For $d = 3$, by part 2 of Theorem 2.4, $V$ hits the set $\Phi^{-1}(E_3)$ with positive probability, since by the absolute continuity result proven above, the hitting properties of $V$ and $U$ over $[\epsilon, T] \times [\epsilon, 1-\epsilon]$ are the same. Therefore, by (6.4), for $\delta = 3$, $u^{(3)}$ hits $E_3$ with positive probability.

Setting

$$E_1 := \{\bar{u} \in C^+ : \exists\, 0 < x < 1,\ \bar{u}(x) = 0\}$$

for $d = 5$, $V$ hits $\Phi^{-1}(E_1)$ with positive probability by Theorem 2.3, so that, for $\delta = 5$ and, by monotonicity, for all $\delta \in [3,5]$, $u^{(\delta)}$ hits $E_1$ with positive probability. $\square$

R. C. DALANG
INSTITUT DE MATHÉMATIQUES
ÉCOLE POLYTECHNIQUE FÉDÉRALE
STATION 8
1015 LAUSANNE
SWITZERLAND
E-MAIL: robert.dalang@epfl.ch

C. MUELLER
DEPARTMENT OF MATHEMATICS
UNIVERSITY OF ROCHESTER
ROCHESTER, NEW YORK 14627
USA
E-MAIL: cmlr@math.rochester.edu

L. ZAMBOTTI
DIPARTIMENTO DI MATEMATICA
POLITECNICO DI MILANO
PIAZZA LEONARDO DA VINCI 32
I-20133 MILANO
ITALY
E-MAIL: zambotti@mate.polimi.it